\documentstyle[12pt]{article}
\newcommand{\sect}[1]{\section{#1}\setcounter{equation}{0}}

\font\mbn=msbm10 scaled \magstep1
\font\mbs=msbm7 scaled \magstep1
\font\mbss=msbm5 scaled \magstep1
\newfam\mbff
\textfont\mbff=\mbn
\scriptfont\mbff=\mbs
\scriptscriptfont\mbff=\mbss\def\mbf{\fam\mbff}

\def\Z{{\mbf Z}}
\def\Co{{\mbf C}}
\def\To{{\mbf T}}
\def\P{{\mbf P}}
\def\Di{{\mbf D}}
\newtheorem{Th}{Theorem}[section]
\newtheorem{Lm}[Th]{Lemma}
\newtheorem{C}[Th]{Corollary}
\newtheorem{D}[Th]{Definition}
\newtheorem{Proposition}[Th]{Proposition}
\newtheorem{R}[Th]{Remark}
\newtheorem{E}[Th]{Example}
\author{Alexander Brudnyi\thanks{Research supported in part by NSERC.
\newline
1991 {\em Mathematics Subject Classification}. Primary 30D15,
Secondary 32F05.
\newline
{\em Key words and phrases}.
Holomorphic vector bundle, Grauert theorem, cohomology, Stein manifold.
}\\
Department of Mathematics and Statistics\\
University of Calgary, Calgary\\
Canada}
\title{A Grauert Type Theorem and Extension of Matrices with Entries
in $H^{\infty}$}
\date{}
\begin{document}
\maketitle
\begin{abstract}
{In the paper we prove an extension theorem for matrices with entries
in $H^{\infty}(U)$ for $U$ being a Riemann surface of a special type. One
of the main components of the proof is a Grauert type theorem for 
``holomorphic'' vector bundles defined over maximal ideal spaces of certain 
Banach algebras.}
\end{abstract}
\sect{\hspace*{-1em}. Introduction.}
{\bf 1.1.} Let $N\subset\subset M$ be a relatively compact domain in an
open Riemann surface $M$ such that
\begin{equation}\label{isom}
\pi_{1}(N)\cong\pi_{1}(M)\ .
\end{equation}
Let $R$ be an unbranched covering of $N$ and $i:U\hookrightarrow R$ be a
domain in $R$. Assume that
\begin{equation}\label{ker}
\begin{array}{l}
\mbox{{\em the induced homomorphism 
of the fundamental groups}}\\
\mbox{{\em $i_{*}:\pi_{1}(U)\longrightarrow\pi_{1}(R)$ is injective.}}
\end{array}
\end{equation}

In this paper we continue to study the space $H^{\infty}(U)$,
of bounded holomorphic functions on $U$ satisfying
(\ref{isom}) and (\ref{ker}), started in [Br]. One of the main results
proved in [Br,Th.1.1] was a Forelli type theorem on projections
for $H^{\infty}(U)$. In the present paper we prove
an extension theorem for matrices with entries in $H^{\infty}(U)$. To 
formulate the result let us recall the following definition.

We say that a collection $f_{1},...,f_{n}$ of functions from
$H^{\infty}(U)$ satisfies the {\em corona condition} if
\begin{equation}\label{coro}
|f_{1}(z)|+|f_{2}(z)|+...+|f_{n}(z)|\geq\delta>0\ \ \ {\rm for\ all}
\ \ \ z\in U\ .
\end{equation}
In [Br,Corol.1.5] we proved that the {\em corona problem} is solvable in
$H^{\infty}(U)$ meaning
that for any $f_{1},...,f_{n}$ satisfying
(\ref{coro}) there are $g_{1},...,g_{n}\in H^{\infty}(U)$ such that
\begin{equation}\label{solu}
f_{1}g_{1}+f_{2}g_{2}+...+f_{n}g_{n}\equiv 1\ .
\end{equation}
For instance, for $U$ being the open unit disc $\Di\subset\Co$ the
solvability of the corona problem follows from the celebrated Carleson's 
Corona Theorem [C]. In this paper we consider a
matrix version of the corona problem.
\begin{Th}\label{te1}
Let $A=(a_{ij})$ be a $n\times k$ matrix, $k<n$, with entries in
$H^{\infty}(U)$. Assume that the family of determinants of submatrices of $A$
of order $k$ satisfies the corona condition. Then there exists an
$n\times n$ matrix $\tilde A=(\tilde a_{ij})$,
$\tilde a_{ij}\in H^{\infty}(U)$, so that
$\tilde a_{ij}=a_{ij}$ for $1\leq j\leq k$, and $det(\tilde A)=1$.
\end{Th}
In fact, we can estimate the norm of $\tilde A$ in 
terms of the norm of $A$, $\delta$ (from (\ref{coro})), $n$ and $N$. (This
estimate does not depend of the choice of $U$.)
\begin{R}\label{re1}
{\rm For $k=1$ we have a column of functions from 
$H^{\infty}(U)$ satisfying the corona condition. The conclusion of the
theorem in this case is essentially stronger than just the solvability
of equation (\ref{solu}). 

Note that a similar to Theorem \ref{te1} result 
for $H^{\infty}(U)$ with 
$U$ being the interior of a bordered Riemann surface was proved first
by Tolokonnikov  [T,Th.3] (see also this paper for further results and
references concerning the extension problem for matrices with entries in
different function algebras).}
\end{R}
Let us give an example of a Riemann surface $U$ satisfying (\ref{isom}) and
(\ref{ker}).
\begin{E}\label{e1}
{\rm Consider the standard action of the group $\Z+i\Z$ on $\Co$
by shifts. The fundamental domain of the action is the square 
$R:=\{z=x+iy\in\Co\ :\ \max\{|x|,|y|\}\leq 1\}$. By $R_{t}$ we denote the
square similar to $R$ with the length of the side $t$. Let $O$ be the orbit of
$0\in\Co$ with respect to the action of $\Z+i\Z$. 
For any $x\in O$ we will choose some
$t(x)\in [1/2,3/4]$ and consider the square $R(x):=x+R_{t(x)}$ centered at 
$x$. Let $V\subset\Co$ be a simply connected domain satisfying the property:

{\em there is a subset\ $\{x_{i}\}_{i\in I}\subset O$\ such that\
$V\cap(\cup_{x\in O}R(x))=\cup_{i\in I}R(x_{i})$.}\\
We set $U:=V\setminus (\cup_{i\in I}\ R(x_{i}))$. Then $U$ satisfies
the required conditions. In fact, the quotient space $\Co/(\Z+i\Z)$ is 
a torus $\Co\To$. Let $S$ be the image of $R_{1/3}$ in $\Co\To$. Then
$U$ belongs to the covering $C$ of $\Co\To\setminus{S}$ with the covering 
group $\Z+i\Z$. The condition that embedding $U\hookrightarrow C$
induces an injective homomorphism of fundamental groups follows from the
construction of $U$.}
\end{E}
{\bf 1.2.} Two essential components of our proof of Theorem \ref{te1} are
[Br,Th.1.1] (a Forelli type theorem), and 
a new Grauert type theorem formulated in this section.

Let $N\subset\subset M$ be a relatively compact domain of a 
connected Stein manifold $M$. Assume that 
\begin{equation}\label{stein}
\mbox{\em
$\overline{N}$ is
holomorphically convex in $M$, and $\pi_{1}(N)\cong\pi_{1}(M)$.}
\end{equation}
Let ${\cal G}$ be a family of subgroups of $\pi_{1}(M)$.
For any $G\in {\cal G}$ by $p_{G}:M_{G}\longrightarrow M$
we denote the unbranched covering of $M$ corresponding to $G$, that is,
$\pi_{1}(M_{G})=G$. Let $U$ be a domain in $M$ satisfying
$\pi_{1}(U)\cong\pi_{1}(M)$. According to the covering homotopy theorem 
(see e.g. [Hu]), we have $U_{G}:=p_{G}^{-1}(U)\subset M_{G}$ is the covering 
of $U$ corresponding to $G$. In particular, it is valid for $N$. Further,
disjoint union $U_{\cal G}:=\sqcup_{G\in {\cal G}}U_{G}$ is an open subset of
the complex space $M_{\cal G}:=\sqcup_{G\in {\cal G}}M_{G}$. 
By $H^{\infty}(U_{\cal G})$ we denote the Banach algebra
of bounded holomorphic functions on $U_{\cal G}$ equipped with the 
supremum norm. Assume now that $U$ is such that $\overline{N}\subset U$. Let 
$r_{U}:H^{\infty}(U_{\cal G})\longrightarrow H^{\infty}(N_{\cal G})$ be the 
restriction homomorphism. By $H^{\infty}(\overline N_{\cal G})$ we denote
the closure in $H^{\infty}(N_{\cal G})$ of the algebra generated by all
$r_{U}(H^{\infty}(U_{\cal G}))$ with $\overline N\subset U$. 

Let ${\cal M}_{\cal G}(N)$ be the maximal ideal space of 
$H^{\infty}(\overline N_{\cal G})$, that is, the set of all non-trivial
homomorphisms $\phi:H^{\infty}(\overline N_{\cal G})\longrightarrow\Co$
equipped with the weak $*$ topology (which is called
the {\em Gelfand topology}). It is
a compact Hausdorff space. Evaluation at an $x\in N_{\cal G}$ determines an
element of ${\cal M}_{\cal G}(N)$. Hence there is a continuous embedding 
$i:N_{\cal G}\hookrightarrow {\cal M}_{\cal G}(N)$. In what follows we 
regard $N_{\cal G}$ as a subset of ${\cal M}_{\cal G}(N)$. 
Then we prove the corona theorem for $H^{\infty}(\overline N_{\cal G})$.
\begin{Th}\label{corona}
$N_{\cal G}$ is an open everywhere dense subset of ${\cal M}_{\cal G}(N)$.
\end{Th}

Let $U\subset\subset M$ be a relatively compact domain such that 
$\overline{N}\subset U$, and $\overline{U}$ satisfies conditions
(\ref{stein}). Clearly ${\cal M}_{\cal G}(N)\subset {\cal M}_{\cal G}(U)$.
Let $E$ be a continuous vector bundle
over ${\cal M}_{\cal G}(U)$ of complex rank $n$. We say that 
$E$ is {\em holomorphic} if $E|_{U_{\cal G}}$ is holomorphic in the usual
sense. A homomorphism $h: E_{1}\longrightarrow E_{2}$ of holomorphic vector
bundles over  ${\cal M}_{\cal G}(U)$ is said to be {\em holomorphic}
if $h|_{U_{\cal G}}:E_{1}|_{\cal U_{G}}\longrightarrow 
E_{2}|_{\cal U_{G}}$ is a holomorphic map. If, in addition,
$h$ is a homeomorphism we say that $E_{1}$ and $E_{2}$ are 
{\em holomorphically isomorphic}. 
\begin{Th}\label{gra}{\bf (A Grauert Type Theorem)}\
Assume that holomorphic vector bundles $E_{1}$ and $E_{2}$ over some
${\cal M}_{\cal G}(U)$ as above are isomorphic as continuous bundles. 
Then their restrictions to ${\cal M}_{\cal G}(N)$ are holomorphically 
isomorphic.
\end{Th}
{\bf 1.3.} In this section we formulate some corollaries of Theorem \ref{gra}
that will be used in the proof of Theorem \ref{te1}.

Let $U_{n}$ denote the group of unitary $n\times n$ matrices. Assume that
$U$ is an open Riemann surface satisfying the conditions of Theorem
\ref{te1}. Clearly, the universal covering of $U$ is the open unit disk 
$\Di\subset\Co$, and
$\pi_{1}(U)$ acts holomorphically on $\Di$ by M\"{o}bius transformations. 
Further, for a Riemann surface $X$ let us denote by
$||f||_{\infty}$ the norm of $f\in H^{\infty}(X)$. We say 
that a matrix $a=(a_{ij})$
with entries in $H^{\infty}(X)$ is (left/right) {\em invertible} if 
$a^{-1}(z)$ exists 
for each $z\in X$ and $a^{-1}$ has entries in $H^{\infty}(X)$. By
$||a||:=\max_{i,j}||a_{ij}||_{\infty}$ we denote the norm of $a$.
\begin{Th}\label{te2}
Let $\rho:\pi_{1}(U)\longrightarrow U_{n}$ be a homomorphism. There
are a constant $C=C(n,N)>0$ depending only on $n$ and $N$ and an invertible 
$n\times n$ matrix $a=(a_{ij})$,
$a_{ij}\in H^{\infty}(\Di)$, such that
\begin{description}
\item[{\rm (1)}]\ $a(g(z))=a(z)\cdot\rho(g)$\ for any\
$g\in\pi_{1}(U)$, $z\in\Di$;
\item[{\rm (2)}]\ $\max\{||a||,||a^{-1}||\}\leq C$
\end{description}
\end{Th}
Now, by $H_{n}^{\infty}(U)$ we denote the $H^{\infty}(U)$-module 
consisting of the columns $(f_{1},...,f_{n})$, $f_{i}\in H^{\infty}(U)$, 
$i=1,...,n$. Any $H^{\infty}(U)$-invariant subspace of $H_{n}^{\infty}(U)$ 
will be called a {\em submodule}. We say that a submodule 
$M\subset H_{n}^{\infty}(U)$ is closed in the topology of the pointwise 
convergence on $U$ if for any net $\{f_{\alpha}\}\subset M$ that pointwise
converges on $U$ to an $f\in H_{n}^{\infty}(U)$ we have $f\in M$.
As an application of Theorem \ref{te2} we obtain.
\begin{Th}\label{modules}
Let $M\subset H_{n}^{\infty}(U)$ be a submodule closed in the
topology of the pointwise convergence on $U$. Then
for some $k$ the module $M$ can be represented as 
$M=H\cdot H_{k}^{\infty}(U)$, where $H$ is a left invertible $n\times k$
matrix with entries in $H^{\infty}(U)$. Moreover, there is a constant
$c=c(n,N)>0$ depending on $n$ and $N$ only such that
$$
\max\{||H||,||H^{-1}||\}\leq c\ .
$$
\end{Th}
Another possible application of Theorem \ref{te2} is the definition of 
analogs of Blaschke products on $U$. 

Let $r:\Di\longrightarrow U$ be the universal
covering map. By ${\cal Z}(g)$ we denote the divisor of zeros of a non-zero
$g$.
\begin{C}\label{te3}
Let $\{z_{i}\}\subset U$ be a sequence of not necessarily distinct points.
Assume that $r^{-1}(\{z_{j}\})={\cal Z}(f)$ for some $f\in H^{\infty}(\Di)$.
Then there are a positive constant $A=A(N)$ depending on $N$ only and a 
function $h\in H^{\infty}(U)$ such that 
$$
{\cal Z}(h)=\{z_{i}\}\ \ \ {\rm and}\ \ \
\sup_{z\in U}|h(z)|\leq A .
$$
\end{C}
%=============================
\sect{\hspace*{-1em}. Proof of the Corona Theorem for 
${\bf H^{\infty}(\overline{N}_{\cal G})}$.}
In this section we will prove Theorem \ref{corona}.\\
{\bf 2.1.} First, we describe $M_{\cal G}$ as a fibre bundle over $M$
with a discrete fibre. We start with the description of the covering
$p_{G}:M_{G}\longrightarrow M$ corresponding to a group 
$G\in {\cal G}$.

Let ${\cal U}=(U_{i})_{i\in I}$ be an open acyclic cover of $M$ by
sets biholomorphic to open Euclidean balls. For a complex Lie group $S$
by $Z_{\cal O}^{1}({\cal U},S)$ we denote the set of holomorphic $S$-valued
${\cal U}$-cocycles. By definition,
$s=\{s_{ij}\}$, $s_{ij}\in {\cal O}(U_{i}\cap U_{j}, S)$, is an element of
$Z_{\cal O}^{1}({\cal U},S)$ if
$$
s_{ij}s_{jk}=s_{ik}\ \ \  {\rm on}\ \ \  U_{i}\cap U_{j}\cap U_{k}.
$$
Now let $X_{G}:=\pi_{1}(M)/G$ be the
set of cosets of $\pi_{1}(M)$ with respect to the (left) action of $G$ on
$\pi_{1}(M)$ defined by left multiplications. By $[Gq]\in X_{G}$ we
denote the coset containing $q\in\pi_{1}(M)$. Consider the complex
Lie group $H(X_{G})$ of all homeomorphisms of $X_{G}$
(equipped with discrete topology). We define the homomorphism 
$\tau_{G}:\pi_{1}(M)\longrightarrow H(X_{G})$ by the formula:
$$
\tau_{G}([Gq]):=[Gqg^{-1}],\ \ \ q\in\pi_{1}(M) .
$$
Set $Q(G):=\pi_{1}(M)/Ker(\tau_{G})$ and let $[g]_{G}$ be the image of
$g\in\pi_{1}(M)$ in $Q(G)$. Finally by
$\tau_{G}':Q(G)\longrightarrow H(X_{G})$ denote the unique homomorphism
whose pullback to $\pi_{1}(M)$ coincides with $\tau_{G}$. Then 
from the basic facts of the theory of fibre bundles (see e.g. [Hi]) it
follows that

{\em There is a cocycle 
$c=\{c_{ij}\}\in Z_{\cal O}^{1}({\cal U},\pi_{1}(M))$
such that $M_{G}$ is biholomorphic to the quotient space of
$\sqcup_{i\in I}U_{i}\times X_{G}$ by the equivalence relation:
$$
U_{i}\times X_{G}\ni x\times \tau_{G}'([c_{ij}]_{G})(h)\sim x\times h\in
U_{j}\times X_{G}.
$$
Projection $p_{G}:M_{G}\longrightarrow M$ is
defined by the coordinate projections
$U_{i}\times X_{G}\longrightarrow U_{i}$.}

Consider now $M_{\cal G}:=\sqcup_{G\in {\cal G}}M_{G}$.
Let $X_{\cal G}:=\prod_{G\in {\cal G}}X_{G}$ and
$H_{\cal G}:=\prod_{G\in {\cal G}}H_{G}$. 
For any $a\in H_{\cal G}$, $a=\{a_{G}\}_{G\in {\cal G}}$, $a_{G}\in H_{G}$,
we define $a:X_{\cal G}\longrightarrow X_{\cal G}$ by the formula
$$
a(x):=\{a_{G}(x_{G})\}_{G\in {\cal G}},\ \ \ x\in X_{\cal G},\
x=\{x_{G}\}_{G\in {\cal G}},\ x_{G}\in X_{G}.
$$ 
Further, set $\tau_{\cal G}':=\{\tau_{G}'\}_{G\in {\cal G}}$. Then
$\tau_{\cal G}'$ is a homomorphism from
$Q({\cal G}):=\prod_{G\in {\cal G}}Q(G)$ into $H_{\cal G}$. We also
set $[c_{ij}]_{\cal G}:=\{[c_{ij}]_{G}\}_{G\in\cal G}\in
{\cal O}(U_{i}\cap U_{j}, Q(\cal G))$.
Now using the above construction of $M_{G}$ we have that
$M_{\cal G}$ is biholomorphic to the quotient space of
$\sqcup_{i\in I}U_{i}\times X_{\cal G}$ by the equivalence relation:
$$
U_{i}\times X_{\cal G}\ni x\times\tau_{\cal G}'([c_{ij}]_{\cal G})(h)
\sim x\times h\in U_{j}\times X_{\cal G}.
$$
In particular, $M_{\cal G}$ is a bundle over $M$ with fibre $X_{\cal G}$. By
$p_{\cal G}:M_{\cal G}\longrightarrow M$ we denote the corresponding
projection. 

Similarly, for a domain $U\subset M$ satisfying 
$\pi_{1}(U)\cong\pi_{1}(M)$, the complex space
$U_{\cal G}:=\sqcup_{G\in {\cal G}}U_{G}\subset M_{\cal G}$ 
is a bundle over $U$ with discrete fibre $X_{\cal G}$.\\
{\bf 2.2.} As the next step we define a compact Hausdorff space 
$E(\overline N,\beta X_{\cal G})$. Then in the proof of Theorem
\ref{corona}, we will show that this space  is homeomorphic to 
${\cal M}_{\cal G}(N)$.

Let $l^{\infty}(X_{\cal G})$ be the algebra of bounded complex-valued
functions $f$ on the discrete space $X_{\cal G}$ with pointwise multiplication
and norm $||f||=\sup_{x\in X_{\cal G}}|f(x)|$. Let $\beta X_{\cal G}$ be 
the {\em Stone-\v{C}ech compactification} of $X_{\cal G}$, i.e., the
maximal ideal space of $l^{\infty}(X_{\cal G})$ equipped with the Gelfand
topology. Then $X_{\cal G}$ is naturally embedded into $\beta X_{\cal G}$
as an open everywhere dense subset, and the topology on $X_{\cal G}$ induced
by this embedding coincides with the original one, i.e., is discrete.
Every function $f\in l^{\infty}(X_{\cal G})$ has a unique extension
$\hat f\in C(\beta X_{\cal G})$. Further, any homeomorphism 
$h\in H_{\cal G}$ of $X_{\cal G}$ determines an isometric 
isomorphism of Banach algebras 
$h^{*}:l^{\infty}(X_{\cal G})\longrightarrow l^{\infty}(X_{\cal G})$. 
Therefore $h$ can be extended to a homeomorphism 
$\hat h:\beta X_{\cal G}\longrightarrow\beta X_{\cal G}$. This extension
shows that now we can think of $H_{\cal G}$ as a subgroup of the group of
homeomorphisms of $\beta X_{\cal G}$.

We retain the notation of Section 2.1. Let us define
the bundle $E(M,\beta X_{\cal G})$ over $M$ with the fibre
$\beta X_{\cal G}$ as the quotient space
of $\sqcup_{i\in I}U_{i}\times\beta X_{\cal G}$ by the equivalence
relation:
$$
U_{i}\times\beta X_{\cal G}\ni x\times
\tau_{\cal G}'([c_{ij}]_{\cal G})(\xi)\sim x\times\xi
\in U_{j}\times\beta X_{\cal G}.
$$
Let $\hat p_{\cal G}: E(M,\beta X_{\cal G})\longrightarrow M$ be the 
corresponding
projection. Since $U_{i}\times X_{\cal G}$ is an open everywhere
dense subset of $U_{i}\times\beta X_{\cal G}$, the 
definitions of $M_{\cal G}$ and $E(M,\beta X_{\cal G})$ show that
$M_{\cal G}$ is an open everywhere dense subbundle of $E(M,\beta X_{\cal G})$
and $\hat p_{\cal G}|_{M_{\cal G}}=p_{\cal G}$. Similarly we can define the
bundle $E(U,\beta X_{\cal G})\longrightarrow U$ for any domain
$U\subset M$ satisfying $\pi_{1}(U)\cong\pi_{1}(M)$, such that
$U_{\cal G}$ is an open everywhere dense subset of $E(U,\beta X_{\cal G})$.

The proof of the next result repeats the arguments of the proof of
[Br1,Th.2.2].
\begin{Proposition}\label{ext}
Let $U\subset M$ be a domain satisfying $\pi_{1}(U)\cong\pi_{1}(M)$.
For every $h\in H^{\infty}(U_{\cal G})$ there is a unique 
$\hat h\in C(E(U,\beta X_{\cal G}))$ such that $\hat h|_{U_{\cal G}}=h$.
\end{Proposition}
We just recall how to construct $\hat h$. Consider the restriction of $h$ to
the set $U_{i}\times X_{\cal G}$. Then for each $z\in U_{i}$ we extend
the function $h(z,\cdot)\in l^{\infty}(X_{\cal G})$ to 
$\hat h(z,\cdot)\in C(\beta X_{\cal G})$ by continuity. The collection 
of all such extended functions for any $z\in U_{i}$ and $i\in I$ forms
the required function $\hat h$.

Let $E(\overline N,\beta X_{\cal G}):=\hat p_{\cal G}^{-1}(\overline N)$.
Then it is easy to see that 
$N_{\cal G}\subset E(\overline N,\beta X_{\cal G})$ is 
an open everywhere dense subset, and $E(\overline N,\beta X_{\cal G})$ is
a Hausdorff compact. We will prove that $E(\overline N,\beta X_{\cal G})$
is homeomorphic to ${\cal M}_{\cal G}(N)$ which gives us the proof of
Theorem \ref{corona}.\\
{\bf 2.3.} Let $U\subset M$ be a domain satisfying 
$\pi_{1}(U)\cong\pi_{1}(M)$. We study some analytic properties of 
$E(U,\beta X_{\cal G})$. First, it follows from the definition that
the base of the topology on $E(U,\beta X_{\cal G})$ consists of the sets
$S_{O,H}$ homeomorphic to $O\times H$, where $O$ is an open subset of $U$
biholomorphic to an open Euclidean ball and $H$ is a clopen subset of 
$\beta X_{\cal G}$. Also, $S_{O,H}\cap U_{\cal G}$ is an open
everywhere dense subset of $S_{O,H}$.
\begin{D}\label{holo}
A function $f\in C(S_{O,H})$ is said to be holomorphic
if its restriction to $S_{O,H}\cap U_{\cal G}$ is holomorphic in the 
usual sense. For any open set $W\subset U_{\cal G}$ the function 
$f\in C(W)$ is holomorphic if its restriction to each $S_{O,H}\subset W$
is holomorphic. The sheaf
of germs of holomorphic on $E(U,\beta X_{\cal G})$ functions will be 
denoted by $\widehat{\cal O}_{U}$. 

A sheaf ${\cal F}$ of $\widehat{\cal O}_{U}$-modules on 
$E(U,\beta X_{\cal G})$ is called
syzygetic if
for any fibre $F$ of $\hat p_{\cal G}$ 
there is an open neighbourhood of $F$ over which ${\cal F}$ admits a 
finite free resolution
\begin{equation}\label{syzyg}
0\longrightarrow (\widehat{\cal O}_{U})^{n_{k}}\longrightarrow\dots
\longrightarrow (\widehat{\cal O}_{U})^{n_{1}}\longrightarrow {\cal F}
\longrightarrow 0
\end{equation}
of sheaves of $\widehat{\cal O}_{U}$-modules.
\end{D}

In the next result by $H^{*}(E(U,\beta X_{\cal G}),{\cal F})$ we
denote the {\em \v{C}ech cohomology} groups with values in the sheaf 
${\cal F}$.
\begin{Th}\label{vanish}
Let $U\subset M$ be a Stein domain satisfying $\pi_{1}(U)\cong\pi_{1}(M)$.
Let ${\cal F}$ be a syzygetic sheaf on $E(U,\beta X_{\cal G})$. Then
for any $i\geq 1$
$$
H^{i}(E(U,\beta X_{\cal G}),{\cal F})=0\ .
$$
\end{Th}
{\bf Proof.}
First, note that for any open cover  ${\cal V}$ of $E(U,\beta X_{\cal G})$
one can find a refinement ${\cal V'}$ of ${\cal V}$ satisfying the properties:
\\
(a) there is an open countable cover of a finite multiplicity 
$\tilde {\cal V}=(\tilde V_{i})_{i\in I}$ of $U$ by 
sets biholomorphic to complex Euclidean balls such that each element of
${\cal V'}$ is a subset of one of $\hat p_{\cal G}^{-1}(\tilde V_{i})$ of
the form $S_{\tilde V_{i},H}$, where $H$ is a clopen subset of 
$\beta X_{\cal G}$, and $\hat p_{\cal G}(S_{\tilde V_{i},H})=\tilde V_{i}$;\\
(b) for a fixed $i$ the sets $S_{\tilde V_{i}, H}\in {\cal V'}$ form a 
finite open cover of $\hat p_{\cal G}^{-1}(\tilde V_{i})$ such that any
two distinct sets of the cover are non-intersecting;\\
(c) ${\cal F}$ admits a finite free resolution  (\ref{syzyg}) over each
$\hat p_{\cal G}^{-1}(\tilde V_{i})$.\\
Existence of the above cover follows from the fact that the fibre
$\beta X_{\cal G}$ is a compact totally disconnected set and that $U$
has a countable exhaustion by compact subsets. Thus it suffices to
calculate the cohomology groups with respect to the family of open covers
${\cal V'}$.

For a fixed cover ${\cal V'}$ let us consider the presheaf
of sections of ${\cal F}$ defined on all possible intersections of sets from
${\cal V'}$ and its direct image ${\cal F'}:=(\hat p_{\cal G})_{*}(\cal F)$ 
with respect to $\hat p_{\cal G}$. Then ${\cal F'}$ is a presheaf generated
by groups of sections $\Gamma(\hat p_{\cal G}^{-1}(O),\cal F)$ where $O$ is 
intersection of some sets $\tilde V_{i}$ from (a). Now property (b) of
${\cal V'}$ shows that the groups $H^{i}({\cal V'},{\cal F})$ and
$H^{i}(\tilde {\cal V},{\cal F}')$ are isomorphic 
(see also the arguments in [Br,Prop.2.4]). 
Moreover, according to the
construction of $E(U,\beta X_{\cal G})$,  the sheaf generated by 
${\cal F}'$ can be
identified with the sheaf of germs of holomorphic sections of a
holomorphic bundle $B$ over $U$ with the fibre $C(\beta X_{\cal G})$. Here
$B$ is the bundle associated with the right isometric action of 
$Q({\cal G})$ on $C(\beta X_{\cal G})$ (the shift of the argument) 
obtained from the right action of $Q({\cal G})$ on $\beta X_{\cal G}$.
Thus by (c) ${\cal F}'$ admits a finite free resolution
$$
0\longrightarrow ({\cal O}_{U}^{C(\beta X_{\cal G})})^{n_{k}}\longrightarrow
\dots
\longrightarrow ({\cal O}_{U}^{C(\beta X_{\cal G})})^{n_{1}}\longrightarrow 
{\cal F}'\longrightarrow 0
$$
of sheaves of ${\cal O}_{U}^{C(\beta X_{\cal G})}$-modules over each 
$\tilde V_{i}$ (it
is the direct image of the corresponding resolution for ${\cal F}$ over
$\hat p_{\cal G}^{-1}(\tilde V_{i})$). 
Here ${\cal O}_{U}^{C(\beta X_{\cal G})}$ is sheaf of germs of holomorphic
$C(\beta X_{\cal G})$-valued functions. Now since $U$ is Stein and
$\tilde {\cal V}$ consists of Stein manifolds, by Bungart [B,Sect.3] and the
classical Ler\'{e} theorem (about calculating cohomology groups by 
acyclic coverings) we have (for any $i\geq 1$)
$$
H^{i}(U,{\cal F}')=H^{i}(\tilde {\cal V},{\cal F}')=0\ .
$$
This completes the proof of the theorem.\ \ \ \ \ $\Box$

We prove now 
\begin{Lm}\label{separat}
The algebra $H^{\infty}(\overline N_{\cal G})$ separates points of
$E(\overline{N},\beta X_{\cal G})$.
\end{Lm}
{\bf Proof.} Let $F_{1},F_{2}$ be fibres of $\hat p_{\cal G}$ over distinct
points $z_{1},z_{2}\in M$. We will show that
{\em for any $f_{i}\in C(F_{i})$, $i=1,2$,
there is a holomorphic function 
$H\in\widehat{\cal O}_{M}(E(M,\beta X_{\cal G}))$ such that
$H|_{F_{i}}=f_{i}$, $i=1,2$.} Since the restriction of a holomorphic 
function defined on $E(M,\beta X_{\cal G})$ to 
$E(\overline{N},\beta X_{\cal G})$ is a function from 
$H^{\infty}(\overline N_{\cal G})$, this will give the required statement.

Consider the restriction homomorphism of the sheaves 
$r:\widehat{\cal O}_{M}\longrightarrow 
{\cal C}_{F_{1}}\sqcup {\cal C}_{F_{2}}$. Kernel $Ker(r)$ is a syzygetic
sheaf (see e.g. [B,Lm.3.3]). Thus by Theorem \ref{vanish}, 
$H^{1}(E(M,\beta X_{\cal G}),Ker(r))=0$. Then the homomorphism of 
the global sections
$$
r:\widehat{\cal O}_{M}(E(M,\beta X_{\cal G}))\longrightarrow C(F_{1})\sqcup
C(F_{2})\ 
$$
is surjective.\ \ \ \ \ $\Box$
\\
{\bf 2.4.} We are ready to prove Theorem \ref{corona}. It is just a routine
to check that the statement of the theorem is equivalent to the following 
one:

{\em Let $U\subset M$ be a Stein domain such that $\pi_{1}(U)\cong\pi_{1}(M)$
and $\overline{N}\subset U$. Let $f_{1},...,f_{n}$ be a collection
of bounded holomorphic functions defined on $E(U,\beta X_{\cal G})$
satisfying the corona condition (\ref{coro}) at each point of
$E(U,\beta X_{\cal G})$.
Then there are an open domain $V\supset\overline{N}$ in $U$ with 
$\pi_{1}(V)\cong\pi_{1}(U)$ and bounded holomorphic functions
$h_{1},...,h_{n}$ defined on $E(V,\beta X_{\cal G})$ such that
$\sum_{i=1}^{n}h_{i}f_{i}\equiv 1$.}

Consider the homomorphism $t:(\widehat{\cal O}_{U})^{n}\longrightarrow
\widehat{\cal O}_{U}$ defined as
$$
t(s_{1x},...,s_{nx}):=f_{1x}s_{1x}+...+f_{nx}s_{nx}\ \ \ 
x\in E(U,\beta X_{\cal G})\ .
$$
Here $s_{x}$ denotes the germ of section $s$ at the point $x$. Let us check
that $Ker(t)$ is a syzygetic sheaf on $E(U,\beta X_{\cal G})$. 

Let $z\in U$ be a point and $F_{z}:=\hat p_{\cal G}^{-1}(z)\subset
E(U,\beta X_{\cal G})$ be the fibre over $z$. Then there is an open 
neighbourhood $O\subset U$ of $z$ biholomorphic to an open Euclidean ball
such that $\hat p_{\cal G}^{-1}(O)$ is disjoint union of sets
$S_{O,H_{i}}$, $i=1,...,n$, so that
$|f_{i}(v)|\geq \delta/n$ for any $v\in S_{O,H_{i}}$. Here $\delta>0$ is
the constant from the corona condition (\ref{coro}) for $f_{1},...,f_{n}$.
(We also admit that for some $i$ the sets $S_{O,H_{i}}$ can be empty.) 
Existence of such $O$ and $S_{O,H_{i}}$ follows from the continuity of 
$f_{1},...,f_{n}$ and the fact that
$F_{z}$ is compact and totally disconnected. Now let us define the
homomorphism $k: (\widehat{\cal O}_{U})^{n-1}\longrightarrow Ker(t)$ over
$\hat p_{\cal G}^{-1}(O)$ by the formula: if $x\in S_{O,H_{i}}$ then
$$
k(s_{1x},...,s_{n-1 x}):=(s_{1x},...,s_{i-1 x},\sum_{l=1}^{i-1}
(-f_{lx}/f_{i x})s_{lx}+\sum_{l=i}^{n-1}(-f_{l+1 x}/f_{i x})s_{lx},
s_{i x},...,s_{n-1 x})\ .
$$
Since $H_{i}$ are clopen non-intersecting subsets of $\beta X_{\cal G}$, we 
can glue together the holomorphic
matrices that define $k$ on each $S_{O,H_{i}}$ to obtain a global 
holomorphic homomorphism over $\hat p_{\cal G}^{-1}(O)$. Clearly,
this homomorphism determines an isomorphism between 
$(\widehat{\cal O}_{U})^{n-1}$ and $Ker(t)$ over
$\hat p_{\cal G}^{-1}(O)$.

Now according to Theorem \ref{vanish},
$H^{1}(E(U,\beta X_{\cal G}),Ker(t))=0$. Since also the homomorphism $t$
is locally surjective (because it is surjective at each point), the standard 
argument (that associates to the short exact sequence of sheaves the long
exact sequence of cohomology groups with values in these sheaves)
shows that the homomorphism of global holomorphic sections 
$(\widehat{\cal O}_{U}(E(U,\beta X_{\cal G})))^{n}\longrightarrow
\widehat{\cal O}_{U}(E(U,\beta X_{\cal G}))$ induced by $t$ is surjective. 
In particular, there are $g_{1},...,g_{n}\in
\widehat{\cal O}_{U}(E(U,\beta X_{\cal G}))$ such that 
$\sum g_{i}f_{i}\equiv 1$. It remains then to restrict $g_{1},...,g_{n}$
to any $E(V,\beta X_{\cal G})$ with $V\supset\overline{N}$,
$\overline{V}\subset U$ and $\pi_{1}(V)\cong\pi_{1}(U)$. 
This restriction produces the required
bounded holomorphic functions $h_{1},...,h_{n}$. 

The above arguments together with Lemma \ref{separat} show that 
${\cal M}_{\cal G}(N)$ is homeomorphic to $E(\overline{N},\beta X_{\cal G})$.

The proof of the theorem is complete.\ \ \ \ \ $\Box$
%=============================================
\sect{\hspace*{-1em}. Proof of the Grauert Type Theorem.}
In this section we will prove Theorem \ref{gra}.

Let $E$ be a holomorphic vector bundle defined in an open neighbourhood
$O$ of ${\cal M}_{\cal G}(N)$. Without loss of generality we can consider
$O$ as $E(U,\beta X_{\cal G})$, where $U$ is a Stein domain in $M$ such that
$\overline{N}\subset U$ and $\pi_{1}(M)\cong\pi_{1}(U)$. Since the fibre of
$E(U,\beta X_{\cal G})$ is totally disconnected, the arguments similar
to those used in the proof of Theorem \ref{corona} (see Section 2.4) show
that the sheaf ${\widehat{\cal O}}(E)$ of germs of holomorphic sections of 
$E$ is coherent. (In fact, in an open neighbourhood of a fibre of 
$E(U,\beta X_{\cal G})$ it is holomorphically isomorphic to 
$(\widehat{\cal O}_{U})^{k}$, where $k=rank(E)$.)
By the same reason the sheaf of germs $S(E;F)$ of holomorphic sections of $E$ 
vanishing on a fibre $F$ of $E(U,\beta X_{\cal G})$ is also coherent.
In particular, $H^{1}(E(U,\beta X_{\cal G}), S(E;F))=0$. Now from the
short exact sequence of sheaves
$$
0\longrightarrow S(E;F)\longrightarrow {\widehat{\cal O}}(E)\longrightarrow
{\widehat{\cal O}}(E)|_{F}\longrightarrow 0
$$
by the standard argument involving the long exact sequence of cohomology
groups we obtain
\begin{Proposition}\label{sepsect}
There are global holomorphic sections $s_{1},...,s_{k}$ of $E$, $k=rank(E)$,
whose restrictions to $F$ give a trivialization of $E|_{F}$.
\ \ \ \ \ $\Box$
\end{Proposition}
Recall that $det(E)$ is a complex rank 1 vector bundle which is determined by
taking the determinant of a cocycle defining $E$.
\begin{Proposition}\label{class}
Assume that $det(E)$ is holomorphically trivial. 
Then there are an open neighbourhood $V\subset O$ of ${\cal M}_{\cal G}(N)$, 
a positive integer $t$, and a
holomorphic map $f_{E}:V\longrightarrow\Co^{t+1}$ such that\\
(1) $f_{E}(V)$ belongs to a closed analytical subset $X_{E}$ of $\Co^{t+1}$ 
defined as the set of zeros of a finite family of holomorphic homogeneous 
polynomials;\\
(2) $X_{E}\setminus\{0\}$ is a smooth manifold;\\
(3) there is a holomorphic vector bundle $E'$ over $X_{E}\setminus\{0\}$
such that $f_{E}^{*}(E')=E$.
\end{Proposition}
{\bf Proof.} From Proposition \ref{sepsect} and
compactness of ${\cal M}_{\cal G}(N)$ it follows that there is an open
neighbourhood $V\subset O$ of ${\cal M}_{\cal G}(N)$ and a
finite number of linearly independent holomorphic sections 
$h_{1},...,h_{n}$ of $E|_{V}$ such that their
restrictions to each point $x\in V$ generate $E_{x}$.
As before  we can choose $V$ in the form 
$E(U,\beta X_{\cal G})$, where $U$ is a Stein domain containing 
$\overline{N}$ such that $\pi_{1}(U)\cong\pi_{1}(M)$.
We now recall the following construction (see, e.g., [GH, Ch.1, Sect. 5]).

Let $W$ be the $n$-dimensional complex vector space of holomorphic sections 
of $E|_{V}$ generated by $h_{1},...,h_{n}$. Since for any 
$x\in V$ sections from $W$ generate the fibre $E_{x}$, the subspace 
$\Lambda_{x}\subset W$ of sections vanishing at $x$ is $n-k$-dimensional.
Below we also use the following definitions.

Let $X$ be a $p$-dimensional complex vector space. By $G(s,X)$ we
denote the corresponding complex Grassmanian, i.e., the set of $s$-dimensional
complex linear subspaces in $X$. We also define the {\em universal bundle}
$S\longrightarrow G(s,X)$ of complex rank $s$, whose fibre over each 
$\Lambda\in G(s,X)$ is the subspace $\Lambda$. Clearly, $S$ is a holomorphic
subbundle of the trivial bundle $G(s,X)\times\Co^{p}$. By $S^{*}$ we
denote the dual bundle which under identification 
$*:G(s,X)\longrightarrow G(p-s,X^{*})$ is isomorphic to the universal bundle
over $G(p-s,X^{*})$.

Now there is a natural map 
$$
\tau_{W}:V\longrightarrow G(n-k,W)=G(k,W^{*})
$$
such that 
$$
E=\tau_{W}^{*}(S^{*})\ \ \ {\rm and}\ \ \ 
W=\tau_{W}^{*}(H^{0}(G(k,n),{\cal O}(S^{*})))\ .
$$
One can express the map $\tau_{W}$ explicitly. Indeed, let $e_{1},...,
e_{k}$ be a local reper of holomorphic sections of $E$ (defined in some 
open subset of $V$). Then for $i=1,...,n$,
$$
h_{i}=\sum_{\alpha=1}^{k} a_{i\alpha}e_{\alpha}
$$
for some holomorphic functions $a_{i\alpha}$.
Under the identification
$G(n-k,W)\cong G(k,W^{*})$, the map $\tau_{W}$ is locally given by
$$
\tau_{W}(x):=\left(
\begin{array}{ccc}
a_{11}(x)&\ldots & a_{1k}(x)\\
\vdots &\ddots &\vdots\\
a_{n1}(x)&\ldots & a_{nk}(x)
\end{array}
\right)\ .
$$
From this formula it is clear that $\tau_{W}$ is holomorphic.
Let $P_{k}:G(k,W^{*})\longrightarrow\P^{t}$,
$t={n\choose k}-1$,
be the Pl\"{u}cker embedding into the projective space. Let 
$\{A_{i}\}_{i\in I}$ be the set of all $k\times k$-minors of the above
matrix $(a_{ij})$. Then locally the map $P_{k}\circ\tau_{W}$ is given
as
$$
(P_{k}\circ\tau_{W})(x):=(A_{1}(x): ... : A_{t+1}(x)) .
$$
This formula shows that $P_{k}\circ\tau_{W}$ is also holomorphic. Moreover,
there is a holomorphic vector bundle $\tilde E$ over the complex projective
manifold $X:=P_{k}(G(k,W^{*}))$ such that 
$(P_{k}\circ\tau_{W})^{*}(\tilde E)=E|_{V}$. Since $det(E)$ is 
holomorphically trivial, the  above functions $A_{s}$ are just the local
representation of global holomorphic functions 
$h_{i_{1}}\wedge...\wedge h_{i_{k}}$ (which are the sections of $det(E)$).
In particular, one can define the holomorphic map 
$f_{E}:V\longrightarrow\Co^{t+1}$, $f_{E}(x)=(A_{1}(x),...,A_{t+1}(x))$,
such that $P_{k}\circ\tau_{W}=
\pi\circ f_{E}$, where $\pi:\Co^{t+1}\longrightarrow\P^{t}$ is the natural
projection. Also the image of $f_{E}$ belongs to the complex manifold 
$X_{k}:=\pi^{-1}(X)\in\Co^{t+1}\setminus\{0\}$. Since $X\subset\P^{t}$ is a 
smooth projective manifold, by Chow's theorem $X_{k}$ is defined as the
set of zeros of a finite family of holomorphic homogeneous polynomials.
Thus $X_{E}:=X_{k}\cup\{0\}\subset \Co^{t+1}$ is the set of zeros of the
same family of polynomials. It remains to set $E':=\pi^{*}(\tilde E)$.
Then according to our construction, $f_{E}^{*}(E')=E|_{V}$.

The proof of the proposition is complete.\ \ \ \ \ $\Box$\\
{\bf Proof of Theorem \ref{gra}.} Let $E_{1}, E_{2}$ be holomorphic
vector bundles of complex rank $k$ defined in an open neighbourhood $O$ of 
${\cal M}_{\cal G}(N)=E(\overline{N},\beta X_{\cal G})$. Assume also that
$E_{1}\cong E_{2}$ as continuous bundles. We will prove that there is an open
neighbourhood $V\subset O$ of ${\cal M}_{\cal G}(N)$ such that $E_{1}|_{V}$
and $E_{2}|_{V}$ are holomorphically isomorphic. 

First, we will prove the
theorem under the additional assumption that the bundles
$det(E_{i})$, $i=1,2$, are holomorphically trivial. 
Then as in Proposition \ref{class} we can construct holomorphic maps 
$\tau_{i}:\overline{U}\longrightarrow G(k,n)\hookrightarrow\P^{t}$ 
(with the same $n$) for some open neighbourhood $U\subset O$ of
${\cal M}_{\cal G}(N)$ such that $E_{i}|_{\overline{U}}=\tau_{i}^{*}(S)$, 
$i=1,2$. Since $\overline{U}$ is a compact and
$E_{1}$, $E_{2}$ are topologically isomorphic,  we can choose $n$ so big 
that $\tau_{1}$ and $\tau_{2}$ are homotopic (see, e.g., Husemoller [H]). 
Denote by $H_{t}:\overline{U}\longrightarrow G(k,n)$, $t\in [0,1]$, this 
homotopy. Then $B:=\{B_{t}:=H_{t}^{*}(S),\ t\in [0,1]\}$, is a continuous
bundle on $\overline{U}\times [0,1]$ such that $B_{0}=E_{1}$ and
$B_{1}=E_{2}$. Let $s_{1},...,s_{n}$ be a basis in
$H^{0}(G(k,n),{\cal O}(S))$. Then $s_{i}(t):=H_{t}^{*}(s_{i})$, $i=1,...,n$,
is a family of linearly independent continuous sections of $B$, such that
for any $t$,  $s_{1}(t),...,s_{n}(t)$ generate each fibre of $B_{t}$.
Moreover, each $H_{t}$ is defined by $s_{1}(t),...,s_{n}(t)$
as in the construction of Proposition \ref{class}. 
Since also $det(B)$ is topologically
trivial (because it determines a continuous homotopy between $det(E_{1})$
and $det(E_{2})$), we obtain a continuous homotopy between 
holomorphic maps  $f_{E_{i}}:\overline{U}\longrightarrow\Co^{t+1}$, $i=1,2$,
that covers $H_{t}$. Indeed, the homotopy map is defined by the family of 
global complex-valued continuous functions
$\{s_{i_{1}}(t)\wedge...\wedge s_{i_{k}}(t)\}$ on $\overline{U}\times[0,1]$.

Now the remaining part of the proof can be extracted from 
Novodvorskii's theorem [No], however, for the sake of completeness we will
present a more detailed exposition.

Consider ${\cal M}_{\cal G}(N)$ as the inverse limit of compacts.
Namely, let $\Gamma$ be the set of all finite collections of holomorphic
functions defined in some open neighbourhoods of $\overline{N}$ such that each
collection from $\Gamma$ contains also the functions from 
$f_{E_{1}}$ and $f_{E_{2}}$. 
Let us fix some order in each $\gamma\in\Gamma$ such that $\gamma$ is
started with the functions from $f_{E_{1}}$ and then from $f_{E_{2}}$. 
For $\gamma_{i}\in\Gamma$, $i=1,2$, we will say that 
$\gamma_{1}\leq\gamma_{2}$ if the ordered set $\gamma_{2}$ contains
$\gamma_{1}$ as an ordered subset. Clearly, each $\gamma\in\Gamma$ can
be considered as a holomorphic map into some $\Co^{l(\gamma)}$, where
$l(\gamma)$ is the cardinality of $\gamma$. Then we set
$N_{\gamma}':=\gamma({\cal M}_{\cal G}(N))\subset\Co^{l(\gamma)}$.
Also by $N_{e}'\in\Co^{2t+2}$ we denote the image of ${\cal M}_{\cal G}(N)$
by the map $e=(f_{E_{1}},f_{E_{2}})$. Further,
by $N_{\gamma}$ we will denote the polynomial hull of $N_{\gamma}'$. 
Clearly, if $\gamma_{1}\leq\gamma_{2}$, the projection 
$\pi_{\gamma_{1}}^{\gamma_{2}}:\Co^{l(\gamma_{2})}\longrightarrow
\Co^{l(\gamma_{1})}$ to the first $l(\gamma_{1})$ coordinates maps 
$N_{\gamma_{2}}$ into $N_{\gamma_{1}}$. According to Proposition
\ref{class} (1), the set $X_{1}:=X_{E_{1}}=X_{E_{2}}\subset\Co^{t+1}$ is 
polynomially convex. Thus $N_{\gamma}$ belongs to the polynomially convex set 
$Y_{\gamma}:=(\pi_{e}^{\gamma})^{-1}(X_{1}\times X_{1})$. 
Let $B^{k}(r)\subset\Co^{k}$ be the open Euclidean ball centered at $0$.
We set
$$
\tilde N_{\gamma}:=Y_{\gamma}\cap(N_{\gamma}+B^{l(\gamma)}(1/l(\gamma)))\ .
$$
Since for $\gamma_{1}\leq\gamma_{2}$ we have $l(\gamma_{1})\leq 
l(\gamma_{2})$, the projection $\pi_{\gamma_{1}}^{\gamma_{2}}$ maps 
$\tilde N_{\gamma_{2}}$ into $\tilde N_{\gamma_{1}}$. Moreover, each 
$\tilde N_{\gamma}$ is an open neighbourhood of $N_{\gamma}$ in $Y_{\gamma}$.
Thus we have the inverse limiting system generated by $\tilde N_{\gamma}$ and
$\pi_{\gamma}^{\beta}$. Since ${\cal M}_{\cal G}(N)$ is the maximal ideal 
space of $H^{\infty}(N_{\cal G})$, the
inverse limit of this system coincides with ${\cal M}_{\cal G}(N)$. 
By $\pi_{\gamma}:{\cal M}_{\cal G}(N)\longrightarrow \tilde N_{\gamma}$ we
denote the corresponding map. Let $h_{1}:\Co^{2t+2}\longrightarrow
\Co^{t+1}$, $h_{2}:\Co^{2t+2}\longrightarrow\Co^{t+1}$ be the projections to 
the first and to the last $t+1$ coordinates, respectively.
Then as we have already proved, $f_{E_{1}}:=h_{1}\circ\pi_{e}$ is homotopic 
to $f_{E_{2}}:=h_{2}\circ\pi_{e}$ inside of the smooth manifold
$X_{1}\setminus\{0\}$. Now it is easy to prove (see, e.g., [L, Lm.1]) that 
there is a $\gamma\in\Gamma$ such that 
$h_{i}\circ\pi_{e}^{\gamma}$ maps $\tilde N_{\gamma}$ into
$X_{1}\setminus\{0\}$, $i=1,2$, and  $h_{1}\circ\pi_{e}^{\gamma}$ is homotopic
to $h_{2}\circ\pi_{e}^{\gamma}$ inside of $X_{1}\setminus\{0\}$. 
In particular, we obtain that $\tilde N_{\gamma}$ belongs to the smooth part 
of $Y_{\gamma}$. Let $F_{i}:=(h_{i}\circ\pi_{e}^{\gamma})^{*}(E')$, $i=1,2$, 
with $E'$ as in Proposition \ref{class} (3). Then $F_{1}$ and $F_{2}$ are
isomorphic as continuous vector bundles and 
$E_{i}=(h_{i}\circ\pi_{e})^{*}(E')$, $i=1,2$. Since
$\tilde N_{\gamma}\subset Y_{\gamma}$ is an open smooth neighbourhood of the 
Stein compact $N_{\gamma}$, there is an open smooth Stein
neighbourhood $X\subset \tilde N_{\gamma}$ of $N_{\gamma}$.
Now $F_{1}|_{X}\cong F_{2}|_{X}$ as continuous bundles and therefore 
by Grauert's theorem [Gr], they are holomorphically isomorphic too. Then
$E_{1}$ and $E_{2}$ are holomorphically isomorphic because they are pullbacks
of $F_{1}|_{X}$ and $F_{2}|_{X}$ by a holomorphic map. 

This completes the proof of the theorem in the case when $det(E_{i})$, 
$i=1,2$, are holomorphically trivial.

Consider now the case when $det(E_{1})$ and $det(E_{2})$ are not 
necessarily holomorphically trivial. From the conditions of the theorem
it follows that $det(E_{1})$ and $det(E_{2})$ are topologically isomorphic.
Thus the holomorphic complex rank 1 vector bundle 
$det(E_{1})\otimes det(E_{2})^{-1}$ is topologically trivial. Now
according to the first part of the proof of the theorem, 
$det(E_{1})\otimes det(E_{2})^{-1}$ is holomorphically trivial too.
Consider now the holomorphic vector bundles 
$H_{i}:=E_{i}\otimes det(E_{i})^{-1}$, $i=1,2$. Then $H_{1}$ and $H_{2}$
are topologically isomorphic and $det(H_{i})$, $i=1,2$,
are topologically trivial.
Again from the first part of the proof it follows that $H_{1}$ and 
$H_{2}$ are holomorphically isomorphic. Since as it was shown above
$det(E_{1})$ and $det(E_{2})$ are holomorphically isomorphic,
the bundles $E_{1}=H_{1}\otimes det(E_{1})$ and 
$E_{2}=H_{2}\otimes det(E_{2})$ are holomorphically isomorphic too.

The proof of the theorem is complete.\ \ \ \ \ $\Box$
%=========================
\sect{\hspace*{-1em}. Some Additional Properties of 
${\bf E(M,\beta X_{\cal G})}$.}
In this section we describe some topological properties of the space
$E(M,\beta X_{\cal G})$ which will be used in the proof of Theorem
\ref{te1}.

Let $H$ be a connected Hausdorff space. 
Assume also that $H$ is locally arcwise connected,
locally simply connected and admits an 
exhaustion by at most countable number of compact subsets. Then the 
fundamental group $\pi_{1}(H)$ is well defined. Now
the bundle $E(H,\beta X_{\cal G})$ over $H$ can be defined using the coverings
corresponding to some family of subgroups $G\in {\cal G}$ exactly 
as in Section 2.2 (with $H$ instead of $M$).
Assume, in addition, that the topological (covering) dimension $dim(H)$
of $H$ is finite.
\begin{Proposition}\label{dimen}
$E(H,\beta X_{\cal G})$ is a paracompact space satisfying
$$
dim(E(H,\beta X_{\cal G}))=dim(H)\ .
$$ 
\end{Proposition}
{\bf Proof.} The above conditions imply that $H$
can be covered by at most countable number of simply connected
relatively compact
subsets $U_{i}$. By the definition over each $U_{i}$ the bundle
$E(H,\beta X_{\cal G})$ is homeomorphic to $U_{i}\times\beta X_{\cal G}$.
Therefore $E(H,\beta X_{\cal G})$ can be covered by at most countable number
of open relatively compact subsets $E(H,\beta X_{\cal G})|_{U_{i}}$ which
implies that $E(H,\beta X_{\cal G})$ is a paracompact space.
Since $dim(U_{i})\leq dim(H)$ and $dim(\beta X_{\cal G})=0$, we obtain
that $dim(U_{i}\times\beta X_{\cal G})\leq dim(H)$. But for any
$x\in \beta X_{\cal G}$ we have $U_{i}\times\{x\}\subset 
U_{i}\times\beta X_{\cal G}$. Thus $dim(U_{i})\leq 
dim(U_{i}\times\beta X_{\cal G})$. Now the required result follows from 
existence of the open cover of $E(H,\beta X_{\cal G})$ by sets 
homeomorphic to $U_{i}\times\beta X_{\cal G}$ and the fact that there is 
some $i$ for which $dim(U_{i})=dim(H)$.\ \ \ \ \ $\Box$

Let $K$ be a connected finite-dimensional locally compact Hausdorff space 
satisfying the same properties as $H$. Assume also that $K$ is homotopically 
equivalent to $H$. Then $\pi_{1}(K)\cong\pi_{1}(H)$ and so 
$E(K,\beta X_{\cal G})$ is well defined too.
\begin{Proposition}\label{homotop}
Under the above assumptions, $E(H,\beta X_{\cal G})$ and 
$E(K,\beta X_{\cal G})$ are homotopically equivalent.
\end{Proposition}
{\bf Proof.} Let $F:H\longrightarrow K$ and
$G:K\longrightarrow H$ be such that $F\circ G:K\longrightarrow K$ and
$G\circ F:H\longrightarrow H$ are homotopic to the identity maps. Since
$F_{*}$ and $G_{*}$ induce isomorphisms of fundamental groups, the covering
homotopy theorem (see e.g. [Hu]) implies that there are
continuous maps $\tilde F: H_{\cal G}\longrightarrow K_{\cal G}$ and
$\tilde G:K_{\cal G}\longrightarrow H_{\cal G}$ that cover $F$ and $G$, 
respectively, such that $\tilde F\circ\tilde G$,
and $\tilde G\circ\tilde F$ are homotopic to the identity maps. (Here
$H_{\cal G}$ and $K_{\cal G}$ are defined similarly to $M_{\cal G}$ from
Section 2.1.)
Locally the maps $\tilde F$ and $\tilde G$ can be described as follows.
Let $U\subset H$ and $V\subset K$ be open simply connected subsets such
that $F(U)\subset V$. Then by the definition $H_{\cal G}$ is homeomorphic to
$U\times X_{\cal G}$ over $U$ and $K_{\cal G}$ is 
homeomorphic to $V\times X_{\cal G}$ over $V$. Moreover, in appropriate
local coordinates we have $\tilde F(u\times x):=F(u)\times x
\in K_{\cal G}|_{V}$ for
$u\times x\in U\times X_{\cal G}\cong H_{\cal G}|_{U}$. A similar description
is valid for $\tilde G$. 
Since the local description is clearly equivariant with
respect to the equivalence relations defining $E(H,\beta X_{\cal G})$ and
$E(K,\beta X_{\cal G})$, we can extend $\tilde F$ and $\tilde G$ to 
continuous maps $\tilde F':E(H,\beta X_{\cal G})\longrightarrow
E(K,\beta X_{\cal G})$ and $\tilde G':E(K,\beta X_{\cal G})\longrightarrow
E(H,\beta X_{\cal G})$. For instance, for $\tilde F'$ the required
extension over $U$ has the from
$$
\tilde F'(u\times\xi):=F(u)\times\xi\in E(K,\beta X_{\cal G})|_{V},\ \ \ 
u\times\xi\in U\times\beta X_{\cal G}\cong E(H,\beta X_{\cal G})|_{U}\ .
$$
It remains to check that $\tilde F'$ and $\tilde G'$ determine a 
homotopic equivalence of $E(H,\beta X_{\cal G})$ and $E(K,\beta X_{\cal G})$.
Consider $G\circ F$. Let $R:H\times [0,1]\longrightarrow H$ be
a continuous map such that $R(.,0)=G\circ F$ and $R(.,1)=id$. Let us 
prove that we can lift $R$ to a homotopy of $E(H,\beta X_{\cal G})$.
As before, by the covering homotopy theorem, there is a homotopy
$\tilde R:H_{\cal G}\times [0,1]\longrightarrow H_{\cal G}$ that covers
$R$. Let $U_{i}\subset H$, $i=1,2$, be open simply
connected sets such that $R(U_{1}\times I)\subset U_{2}$ for an open
subinterval $I\subset [0,1]$. Then in appropriate local
coordinates $\tilde R:H_{\cal G}|_{U_{1}}\times I\longrightarrow 
H_{\cal G}|_{U_{2}}$ is defined by the formula
$$
\tilde R(u_{1}\times x\times t):=R(u_{1}\times t)\times x\in 
H_{\cal G}|_{U_{2}},\ \ \ \ u_{1}\times x\times t\in 
H_{\cal G}|_{U_{1}}\times I\ .
$$ 
Clearly, $\tilde R(.,0)=\tilde G\circ\tilde F$ and $\tilde R(.,1)=id$.
Since $\tilde R$ is equivariant with respect to the equivalence
relation that defines $E(H,\beta X_{\cal G})$, we can determine the homotopy
$\tilde R':E(H,\beta X_{\cal G})\longrightarrow E(H,\beta X_{\cal G})$.
Locally it is given by the formula
$$
\tilde R'(u_{1}\times\xi\times t):=R(u_{1}\times t)\times\xi\in 
E(H,\beta X_{\cal G})|_{U_{2}},\ \ \ \ u_{1}\times\xi\times t\in 
E(H,\beta X_{\cal G})|_{U_{1}}\times I\ .
$$ 
Clearly $\tilde R'$ is a continuous extension of $\tilde R$, and
$\tilde R'(.,0)=\tilde G'\circ\tilde F'$, \ $\tilde R'(.,1)=id$. Similar
arguments can be  applied to $\tilde F'\circ\tilde G'$. This shows that
$E(H,\beta X_{\cal G})$ and $E(K,\beta X_{\cal G})$ are homotopically
equivalent.\ \ \ \ \  $\Box$

In the next result we will assume that $M$ is a complex manifold.
Recall that $M_{\cal G}$ is a fibre bundle over $M$ with the fibre
$X_{\cal G}$ associated with a representation $\tau_{\cal G}'$ of
$Q(\cal G)$ into the group of bijective homomorphisms of $X_{\cal G}$.
Let $T:=\prod_{G\in {\cal G}}G$ and
$\rho: T\longrightarrow U_{n}$ be a unitary representation.
By $E_{\rho}$ we denote the flat vector bundle over
$M_{\cal G}$ associated with $\rho$.
\begin{Proposition}\label{bunext}
There is a unique holomorphic vector bundle $\tilde E_{\rho}$ over
$E(M,\beta X_{\cal G})$ whose restriction to $M_{\cal G}$ coincides with
$E_{\rho}$.
\end{Proposition}
{\bf Proof.} 
According to the construction of Section 2.1, we can describe
$E_{\rho}$ as follows (cf. [Br,Prop.2.3]).
There is an open acyclic cover ${\cal U}=(U_{i})_{i\in I}$ of $M$ such that 
$E_{\rho}$ is defined on the cover ${\cal V}=(V_{i})_{i\in I}$ of 
$M_{\cal G}$, $V_{i}:=p_{\cal G}^{-1}(U_{i})$, by a cocycle 
$d=\{d_{ij}\}\in Z_{\cal O}^{1}({\cal V}, T)$, i.e., $E_{\rho}$ is defined by 
the equivalence relation
$$
V_{i}\times\Co^{n}\ni x\times\rho(d_{ij})(v)\sim x\times v\in V_{j}\times
\Co^{n}\ .
$$
Observe that each $V_{i}$ is biholomorphic to $U_{i}\times X_{\cal G}$
and therefore $\rho(d_{ij})$ can be thought of as a function defined on
$U_{i}\cap U_{j}$ with values in the space of maps 
$X_{\cal G}\longrightarrow U_{n}$. In fact, $\rho(d_{ij})$ is
a constant multivalued function. Further, because $U_{n}$ is a compact
subset of some $\Co^{N}$, each map 
$r: X_{\cal G}\longrightarrow U_{n}$ admits the natural continuous extension
$\tilde r: \beta X_{\cal G}\longrightarrow U_{n}$ (obtained by the
extension of coordinates of $r$). Therefore we obtain an extended function
$\tilde\rho(d_{ij}):U_{i}\cap U_{j}\longrightarrow
C(\beta X_{\cal G}, U_{n})$. Now denote by $\overline{V}_{k}$ the set
obtained from $V_{k}$ by taking the closure of each fibre of 
$V_{k}\subset E(M,\beta X_{\cal G})$. By definition, 
$\overline{V}_{k}=\hat{p}_{\cal G}^{-1}(U_{k})\cong  
U_{k}\times\beta X_{\cal G}$ is an open subset of $E(M,\beta X_{\cal G})$.
Rewriting according to this identification 
$\tilde\rho(d_{ij})$ in coordinates on $E(M,\beta X_{\cal G})$, we can
think of $\tilde\rho(d_{ij})$ as a continuous function 
$\overline{V}_{i}\cap\overline{V}_{j}\longrightarrow U_{n}$ such that
$\tilde\rho(d_{ij})|_{V_{i}\cap V_{j}}=\rho(d_{ij})$.
Moreover, since
$\rho(d_{ij})\cdot\rho(d_{jk})=\rho(d_{ik})$ on $V_{i}\cap V_{j}\cap V_{k}$,
we have $\tilde\rho(d_{ij})\cdot\tilde\rho(d_{jk})=\tilde\rho(d_{ik})$
on $\overline{V}_{i}\cap\overline{V}_{j}\cap\overline{V}_{k}$
by continuity. This shows that $\{\tilde\rho(d_{ij})\}$
determines a cocycle defined on the cover $(\overline{V}_{i})_{i\in I}$ of
$E(M,\beta X_{\cal G})$ with values in $U_{n}$. 
Thus we can define $\tilde E_{\rho}$ on
$E(M,\beta X_{\cal G})$ by the equivalence relation
$$
\overline{V}_{i}\times\Co^{n}\ni x\times\tilde\rho(d_{ij})(v)\sim x\times v
\in\overline{V}_{j}\times\Co^{n}\ .
$$
Clearly $\tilde E_{\rho}|_{M_{\cal G}}=E_{\rho}$ and so $\tilde E_{\rho}$ is
holomorphic by the definition.\ \ \ \ \ $\Box$

Assume now that $M$ is a non-compact complex Riemann surface and $E_{\rho}$,
$\tilde E_{\rho}$ are the bundles from Proposition \ref{bunext}.
\begin{Proposition}\label{triviality}
$\tilde E_{\rho}\longrightarrow E(M,\beta X_{\cal G})$ is a topologically 
trivial holomorphic vector bundle.
\end{Proposition}
{\bf Proof.} Since any non-compact complex Riemann surface is a Stein 
manifold,
$M$ is homotopically equivalent to a one-dimensional CW-complex $K$ (which 
satisfies conditions of Proposition \ref{homotop}). Then according
to Proposition \ref{homotop}, the paracompact spaces
$E(M,\beta X_{\cal G})$ and 
$E(K,\beta X_{\cal G})$ are homotopically equivalent. Let 
$s:E(K,\beta X_{\cal G})\longrightarrow E(M,\beta X_{\cal G})$ be one of the
maps determining the equivalence. Then according to the general theory
of vector bundles (see [H]) the statement of the proposition will
follow from triviality of $s^{*}(\tilde E_{\rho})$ over
$E(K,\beta X_{\cal G})$. Now  Proposition \ref{dimen} implies that
$dim(E(K,\beta X_{\cal G}))=1$. Therefore
the only obstacle to triviality of $s^{*}(\tilde E_{\rho})$ is
the first Chern class 
$c_{1}(s^{*}(\tilde E_{\rho}))\in H^{2}(E(K,\beta X_{\cal G}),\Z)$. Since
$E(K,\beta X_{\cal G})$ is one dimensional, the latter cohomology group
is trivial which implies that $c_{1}(s^{*}(\tilde E_{\rho}))=0$, and so
$\tilde E_{\rho}$ is topologically trivial.\ \ \ \ \  $\Box$
%========================
\sect{\hspace*{-1em}. Proof of Theorem 1.4 and Corollary 1.6.}
{\bf Proof of Theorem \ref{te2}.} 
Let $N\subset\subset M$ be a relatively compact domain in an
open Riemann surface $M$ such that
$\pi_{1}(N)\cong\pi_{1}(M)$.
Let $R$ be an unbranched covering of $N$ and $i:U\hookrightarrow R$ be a
domain in $R$. Assume that the induced homomorphism 
of the fundamental groups
$i_{*}:\pi_{1}(U)\longrightarrow\pi_{1}(R)$ is injective. Without loss
of generality we may also assume that $i_{*}$ is surjective. Indeed, if
$i_{*}(\pi_{1}(U))$ is a proper subgroup of $\pi_{1}(R)$, then consider
the covering $\tilde R$ of $R$ corresponding to $i_{*}(\pi_{1}(U))$.
Clearly $\tilde R$ is a covering of $N$. Now, by the covering homotopy
theorem, there is a holomorphic embedding 
$\tilde i:U\hookrightarrow\tilde R$ that covers $i$ and such that
$\tilde i_{*}:\pi_{1}(U)\longrightarrow\pi_{1}(\tilde R)$ is a bijection.
So we can work with the triple $(U,\tilde i,\tilde R)$ satisfying conditions
of the theorem.

Now, according to the above assumptions, any homomorphism
$\rho:\pi_{1}(R)\longrightarrow U_{n}$ coincides with
$\rho\circ i_{*}:\pi_{1}(U)\longrightarrow U_{n}$. 

Let $G\subset\pi_{1}(N)$ be a subgroup and $\rho: G\longrightarrow U_{n}$ be
a homomorphism. Denoting $G$ by $G_{\rho}$ we emphasize that $G$ is 
the domain of the definition of $\rho$.
Then we set 
$$
{\cal G}:=\{G_{\rho}\ :\ G\subset\pi_{1}(N),\ \rho\in Hom(G,U_{n})\}\ .
$$
Since $\pi_{1}(M)\cong\pi_{1}(N)$, we can construct
the associated with a $\rho\in Hom(G,U_{n})$ complex vector bundle 
$E_{\rho}$ over $M_{G_{\rho}}$. Here $M_{G_{\rho}}$ is the covering of 
$M$ with the fundamental group $G_{\rho}$. 
Let ${\cal U}=(U_{i})_{i\in I}$ be the acyclic cover 
of $M$ from the 
construction of $M_{G_{\rho}}$ and $M_{\cal G}$ (see Section 2.1). 
Consider the cover ${\cal V}_{G_{\rho}}:=p_{G_{\rho}}^{-1}({\cal U})$ of 
$M_{G_{\rho}}$.
Since ${\cal V}_{G_{\rho}}$ is acyclic too, there is a cocycle
$\{c_{ij, G_{\rho}}\}\in Z_{\cal O}^{1}({\cal V}_{G_{\rho}}, U_{n})$ such that
$E_{\rho}$ is equivalent to the quotient space of 
$\sqcup_{i\in I}\ p_{G_{\rho}}^{-1}(U_{i})\times\Co^{n}$ by the equivalence 
relation
$$
p_{G_{\rho}}^{-1}(U_{i})\times\Co^{n}\ni x\times c_{ij,G_{\rho}}(v)\sim x
\times v\in p_{G_{\rho}}^{-1}(U_{j})\times\Co^{n}\ .
$$
The same is valid for any $\rho$ and $G_{\rho}$. Thus we can construct a 
holomorphic 
bundle $E_{\cal G}$ over $M_{\cal G}$ defined on the cover 
$p_{\cal G}^{-1}({\cal U})$ such that $E_{\cal G}|_{M_{G_{\rho}}}=E_{\rho}$.
Clearly, the bundle $E_{\cal G}$ has the unitary structure group. In fact,
$E_{\cal G}$ is a vector bundle associated with a representation
$R:\prod_{G_{\rho}\in {\cal G}}G_{\rho}\longrightarrow U_{n}$, 
$R|_{G_{\rho}}=\rho$.
Now according to Propositions \ref{bunext} and \ref{triviality}, 
the extension 
$\tilde E_{\cal G}$ of $E_{\cal G}$ to $E(M,\beta X_{\cal G})$ is 
topologically trivial. 
Since $\overline{N}\subset M$ is a Stein compact, by Theorem \ref{gra}
$\tilde E_{\cal G}$ is holomorphically trivial in an open neighbourhood 
of $E(\overline{N},\beta X_{\cal G})\subset E(M,\beta X_{\cal G})$. 
Going back to $E_{\rho}$ we obtain that
cocycle $\{c_{ij,G_{\rho}}\}$ is holomorphically trivial on 
$\overline{N}_{G_{\rho}}\subset M_{G_{\rho}}$.
We can say even more. Let ${\cal U}'=(U_{i}')$ be a finite acyclic
cover of $\overline{N}$ by compact Euclidean balls such that each 
$U_{i}'$ is a compact subset of one of $U_{j}\in {\cal U}$. Let us
restrict $\{c_{ij,G}\}$ to $p_{G_{\rho}}^{-1}({\cal U}')$. From
the triviality of $\tilde E_{\cal G}$ it follows that there are holomorphic
functions $c_{i,G_{\rho}}\in {\cal O}(p_{G_{\rho}}^{-1}(U_{i}'),GL_{n}(\Co))$
such that
$$
c_{i,G_{\rho}}^{-1}(z)\cdot c_{j,G_{\rho}}(z)=c_{ij,G_{\rho}}(z),\ \ \ \ 
{\rm for\ any}\ \ z\in p_{G_{\rho}}^{-1}(U_{i}')\cap 
p_{G_{\rho}}^{-1}(U_{j}'),
$$
and there is a constant $C>0$ depending only on $n$ and $N$ such that
$$
\sup_{i,G_{\rho}}\max\{||c_{i,G_{\rho}}||,||c_{i,G_{\rho}}^{-1}||\}\leq C\ .
$$
(Here $||.||$ is defined as in Section 1.3.)

Let us take now an unbranched covering $R$ of $N$ and $U\subset R$
as in the beginning of this section. 
Consider the universal covering $r:\Di\longrightarrow R$.
(Recall that $R\subset M_{G_{\rho}}$ and $\pi_{1}(R)=
\pi_{1}(M_{G_{\rho}})=G_{\rho}$.)
Then $r^{-1}(p_{G_{\rho}}^{-1}(U_{i}'))=\sqcup_{g\in G_{\rho}}S_{ig}$, where
$S_{ig}$ is biholomorphic to $p_{G_{\rho}}^{-1}(U_{i}')$. 
Consider the pullback
$r^{*}(\{c_{ij,G_{\rho}}\})$ to the cover $(S_{ig})$ of $\Di$. Since the 
latter is an acyclic cover and the bundle $E_{\rho}$ is obtained from the
representation $\rho:G_{\rho}\longrightarrow U_{n}$, there are locally 
constant 
functions $\rho_{ig,G_{\rho}}:S_{ig}\longrightarrow U_{n}$ such that 
$\rho_{ig,G_{\rho}}(z)\cdot\rho_{jh,G_{\rho}}^{-1}(z)=
r^{*}(c_{ij,G_{\rho}})(z)$ for any
$z\in S_{ig}\cap S_{jh}\neq\emptyset$, and $\rho_{il\cdot g,G_{\rho}}=
\rho_{ig,G_{\rho}}\cdot\rho(l)$.
Define $F_{ig,G_{\rho}}=r^{*}(c_{i,G_{\rho}})\cdot\rho_{ig,G_{\rho}}$ on 
$S_{ig}$. 
Then from the above
equations it follows that $F_{ig,G_{\rho}}(z)=F_{jh,G_{\rho}}(z)$ for any 
$z\in S_{ig}\cap S_{jh}$. That is the family $\{F_{ig,G_{\rho}}\}$ 
determines a
global function $F_{G_{\rho}}\in {\cal O}(\Di,GL_{n}(\Co))$. Now the group 
$G_{\rho}$
acts holomorphically on $\Di$ by M\"{o}bius transformations such that
each $l\in G_{\rho}$ maps any $S_{ig}$ biholomorphically onto $S_{il\cdot g}$.
In particular, for $z\in S_{ig}$ we have 
$$
F_{G_{\rho}}(l(z))=r^{*}(c_{i,G_{\rho}})(l(z))\cdot
\rho_{il\cdot g,G_{\rho}}=F_{G_{\rho}}(z)\cdot
(\rho_{ig,G_{\rho}}^{-1}\cdot\rho_{il\cdot g,G_{\rho}})=
F_{G_{\rho}}(z)\cdot\rho(l)\ .
$$
Clearly $F_{G_{\rho}}$ satisfies the required estimates of Theorem \ref{te2}.
Since the universal covering $\tilde U$ of $U$ admits a holomorphic 
embedding into $\Di$ equivariant with respect to the actions of $G_{\rho}$ on
$\tilde U$ and $\Di$, respectively, the restriction $F_{G_{\rho}}|_{\tilde U}$
determines the required matrix function $a$.
\ \ \ \ \ $\Box$
\\
{\bf Proof of Corollary \ref{te3}.} Let $\{z_{i}\}\subset U$ be a sequence
such that $r^{-1}(\{z_{i}\})\subset\Di$ is the set of zeros of a non-zero
bounded holomorphic function. In particular, there is a Blaschke product
$B$ whose set of zeros (counted with their multiplicities) is exactly 
$r^{-1}(\{z_{i}\})$. Since the set $r^{-1}(\{z_{i}\})$ is invariant with
respect to the action of $G:=\pi_{1}(U)$ and $B$ is an interior function, 
we have $B(g(z))=B(z)\cdot\rho(g)$, $g\in G$,
for some representation $\rho:G\longrightarrow U_{1}$. Now according to
Theorem \ref{te2}, there is a holomorphic function 
$a\in {\cal O}(\Di,\Co^{*})$ such that $a(g(z))=a(z)\cdot\rho(g)$, $g\in G$,
and $\max\{||a||_{\infty},||a^{-1}||_{\infty}\}\leq C$ for some 
constant $C$ depending on $N$ only. Set $\tilde h(z):=B(z)/a(z)$. Then 
$\tilde h$ is invariant with respect to the action of $G$, the
set of zeros (counted with their multiplicities) of $\tilde h$ is
$r^{-1}(\{z_{i}\})$, and $||\tilde h||_{\infty}\leq C$. Clearly, 
$\tilde h=r^{*}(h)$ for some $h\in H^{\infty}(U)$ satisfying the 
required properties.\ \ \ \ \ $\Box$
%=====================================
\sect{\hspace*{-1em}. Proof of Theorems 1.1 and 1.5. }
{\bf Proof of Theorem \ref{modules}.}
In the proof we use the Lax-Halmos theorem (see e.g. [T]). By $S^{1}$ we
denote the boundary of $\Di$.\\
{\bf Lax-Halmos Theorem.}
{\em Let $M$ be a weak $*$ closed submodule of the $H^{\infty}(\Di)$-module
$H_{n}^{\infty}(\Di)$. Then for some $k$ we have 
$M=\Psi\cdot H_{k}^{\infty}(\Di)$, where $\Psi$ is a left unimodular
$n\times k$ matrix with entries in $H^{\infty}(\Di)$, that is, 
$\Psi^{*}(\xi)\cdot\Psi(\xi)=I_{k}$\ for a.e. $\xi$, $\xi\in S^{1}$. 
If two such modules 
$\Psi\cdot H_{k}^{\infty}(\Di)$ and $\Theta\cdot H_{m}^{\infty}(\Di)$ are
equal, then $k=m$ and $\Psi=\Theta\cdot V$, where $V\in U_{k}$.}

The proof of a similar result for submodules of $H_{n}^{2}(\Di)$ (with the
same conclusion but with $H_{k}^{2}(\Di)$ instead of $H_{k}^{\infty}(\Di)$)
can be found, e.g., in [Ni,Lect.I, Corol.6]. The required result now can
be obtained from the case of $H_{n}^{2}(\Di)$ similarly to the proof in
[Ga,Ch.II,Th.7.5].

We proceed to the proof of Theorem \ref{modules}. 
In our proof we use a scheme
suggested in [T]. Here, however, instead of the Forelli theorem [F] we use
Theorem 1.1 of [Br], and instead of the classical 
Grauert theorem [Gr] our Theorem \ref{te2}.

Let $r:\Di\longrightarrow U$ be the 
universal covering map
and $G=\pi_{1}(U)$. We can identify $H^{\infty}(U)$ with 
$H_{G}^{\infty}:=r^{*}(H^{\infty}(U))$, the subalgebra 
of $G$-invariant functions in the algebra $H^{\infty}(\Di)$. The module
$M$ can be identified with an $H^{\infty}_{G}$-submodule in 
$H_{G,n}^{\infty}:=r^{*}(H_{n}^{\infty}(U))$. 
Let $N$ be the weak $*$ closed $H^{\infty}(\Di)$-module
generated by $M$. Then by the Lax-Halmos theorem,
$N=\Psi\cdot H_{k}^{\infty}(\Di)$, where $\Psi$ is a left unimodular matrix.
If $A\in G$, then 
$\Psi\cdot H_{k}^{\infty}(\Di)=(\Psi\circ A)\cdot H_{k}^{\infty}(\Di)$ and 
so by the Lax-Halmos theorem $\Psi(A(z))=\Psi(z)\cdot\alpha(A)$ for all
$z\in\Di$, where $\alpha(A)\in U_{k}$. It is clear that
$\alpha:G\longrightarrow U_{k}$ is a homomorphism. According to Theorem
\ref{te2}, there is a bounded matrix $\Omega\in {\cal O}(\Di, GL_{k}(\Co))$
such that $\Omega(A(z))=\Omega(z)\cdot\alpha(A)$. Since
$\alpha(A)^{*}=\alpha(A)^{-1}$, we have
$\Omega'(A(z))=\alpha(A)^{*}\cdot\Omega'(z)$ where by  Theorem
\ref{te2}, $\Omega':=\Omega^{-1}$ is a bounded matrix from
${\cal O}(\Di, GL_{k}(\Co))$. Let 
$$
H^{\infty}(\Di,\alpha):=\{f\in H_{k}^{\infty}(\Di)\ :\
f(A(z))=\alpha(A)^{*}\cdot f(z)\ \ {\rm for\ all}\ z\in\Di\ \ {\rm and\
all}\  A\in G\}\ .
$$
Then $N\cap H_{G,n}^{\infty}=\Psi\cdot H^{\infty}(\Di,\alpha)$.
Further, $\Psi\cdot\Omega'$ is a bounded $n\times k$ matrix with entries in
$H_{G}^{\infty}$, $N\cap H_{G,n}^{\infty}=
\Psi\cdot\Omega'\cdot H_{G,k}^{\infty}$, and $\Psi\cdot\Omega'$
is the pullback of a left invertible $n\times k$ matrix $H$ with
entries in $H^{\infty}(U)$. It remains to check the equality
$N\cap H_{G,n}^{\infty}=M$.\\
By [Br,Th.1.1] there exists a continuous
$H_{G}^{\infty}$-linear projector
$P:H^{\infty}(\Di)\rightarrow H_{G}^{\infty}$ satisfying the property:

{\em Let $o(x):=\{gx\}_{g\in G}$ be an orbit, and let 
$\{f_{\alpha}\}\subset H^{\infty}(\Di)$ be a net.
Assume that the restriction $\{f_{\alpha}\}|_{o(x)}$ converges in 
the weak $*$ topology of $l^{\infty}(o(x))$ to $f|_{o(x)}$ for some
$f\in H^{\infty}(\Di)$. Then $\lim_{\alpha}P(f_{\alpha})(x)=P(f)(x)$. }

Now if $f\in N\cap H_{G,n}^{\infty}$, then $f$ is the limit in the weak $*$
topology of a net $\{f_{\alpha}\}$ of the form 
$f_{\alpha}=\sum_{i=1}^{n_{\alpha}}g_{i\alpha}\cdot h_{i\alpha}$ where
$g_{i\alpha}\in H^{\infty}(\Di)$, $h_{i\alpha}\in M$. Since
any orbit $o(x)$ is an interpolating sequence for $H^{\infty}(\Di)$
(see, e.g. [Br1,Lm.7.1]), the
restriction of $\{f_{\alpha}\}$ to $o(x)$ converges in the weak $*$
topology of $l^{\infty}(o(x))$ to $f|_{o(x)}$. Further,
$P(g_{i\alpha}\cdot h_{i\alpha})=P(g_{i\alpha})\cdot h_{i\alpha}$, and so
we have $P(f_{\alpha})\in M$. Finally, $P(f)=f\in M$ because, according
to our assumptions, $M$ is closed
in the topology of the pointwise convergence on $\Di$. 

To estimate the norms of $H$ and $H^{-1}$ it suffices to estimate the norms
of the matrix $r^{*}(H)=\Psi\cdot\Omega'$ and its inverse. But then the
required estimates follow from the fact that $\Psi$ is a left unimodular 
matrix and $\max\{||\Omega'||,||(\Omega')^{-1}||\}\leq C(n,N)$. 
(Here $||\cdot||$ is defined as in Section 1.3.) Therefore we have 
$\max\{||r^{*}(H)||,||r^{*}(H^{-1})||\}\leq n^{3/2}C(n,N)$.
\ \ \ \ \ $\Box$\\
{\bf Proof of Theorem \ref{te1}.} Let $A=(a_{ij})$ be a 
$n\times k$ matrix, $k<n$, with entries in $H^{\infty}(U)$. Assume
that the family of determinants of submatrices $A$ of order $k$ satisfies
the corona condition. According to Lemma 1 of [T], one can find a
$k\times n$ matrix $G$ with entries in $H^{\infty}(U)$ such that
\begin{equation}\label{ik}
G\cdot A=I_{k}\ .
\end{equation}
The operator $G$ maps $v\in H_{n}^{\infty}(U)$ into
$G\cdot v\in H_{k}^{\infty}(U)$. Let $Ker(G)\subset H_{n}^{\infty}(U)$
be its kernel. Clearly, $Ker(G)$ is a submodule of $H_{n}^{\infty}(U)$
closed in the topology of the pointwise convergence on $U$ (see the definition
in Section 1.3). Then according to Theorem \ref{modules} and (\ref{ik}), 
$Ker(G)=H\cdot H_{n-k}^{\infty}(U)$ for some left invertible 
$n\times (n-k)$ matrix $H$ with entries in $H^{\infty}(U)$. For the matrix
$J:=I_{n}-A\cdot G$, we have $G\cdot J=0$. Let us define $J_{i}$ as the
$i^{\rm th}$ column of the matrix $J$. Then $J_{i}\in H_{n}^{\infty}(U)$,
$G\cdot J_{i}=0$, $J_{i}\in Ker(G)=H\cdot H_{n-k}^{\infty}(U)$, and for
some column $S_{i}\in H_{n-k}^{\infty}(U)$ we have $J_{i}=H\cdot S_{i}$.
Consider the $(n-k)\times n$ matrix $S=(S_{1},...,S_{n})$. Then
$J=H\cdot S$ and $A\cdot G+H\cdot S=I_{n}$. So the $n\times n$ matrix 
$F:=(A,H)$ is invertible and $F^{-1}={G\choose S}$. Since $k<n$, we can
divide the last column of $F$ by $det(F)$ to obtain the invertible
matrix $\tilde A$ with $det(\tilde A)=1$ that extends $A$.

Now, let $\delta>0$ be the number in the corona condition (\ref{coro})
for the minors of $A$ of order $k$ and $||A||$ be the norm of $A$.
Then the effective estimate for solutions of the corona problem in 
$H^{\infty}(U)$ see [Br,Corol.1.5] and Lemma 1 of [T] imply that 
$||G||\leq c_{1}(n,N,\delta,||A||)$. From here and the estimate for
$||H||, ||H^{-1}||$ in Theorem
\ref{modules}, we obtain that $||S||\leq c_{2}(n,N,\delta,||A||)$. 
These estimates together
imply
that $\max\{||\tilde A||,||(\tilde A)^{-1}||\}\leq c_{3}(n,N,\delta,||A||)$
(which does not depend of the choice of $U$).

The proof of the theorem is complete.\ \ \ \ \ $\Box$
%==================================

\end{document}